\documentclass[]{article}
\usepackage{amsmath}
\usepackage{amssymb}
\usepackage{bbm}
\usepackage{hyperref}
\usepackage{tikz}
\usetikzlibrary{decorations.pathmorphing,arrows}
\usepackage{braket}
\usepackage{setspace}

\newtheorem{thm}{Theorem}[section]

\newtheorem{lemma}[thm]{Lemma}

\newtheorem{theorem}[thm]{Theorem}

\newtheorem{corollary}[thm]{Corollary}

\newcommand{\R}{\mathbb{R}}

\newcommand{\Z}{\mathbb{Z}

}

\renewcommand{\Z}{{\mathbb{Z}}}

\title{Online minimum search for a Brownian bridge
}
\author{Erik Wu\footnote{\href{mailto:erikwu35758@gmail.com}{erikwu35758@gmail.com}, }\quad
 and
Shannon Starr\footnote{\href{mailto:slstarr@uab.edu}{slstarr@uab.edu}}
\\
\small
\raisebox{5pt}{$*$}
\large James Clemens High School 
\small \textcolor{white}{$*$}\\
\small 11306 County Line Rd,
Madison, AL 35756\\[3pt]
\small \raisebox{5pt}{$\dagger$}
\large Department of Mathematics, University of Alabama at Birmingham
\small
\textcolor{white}{$\dagger$}\\
\small 1402 10th Avenue South,
Birmingham, AL 35294--1241
}
\date{\today}

\doublespace

\begin{document}

\maketitle

\abstract{In this short note we consider the computational problem of numerically finding the minimum and arg-min of a Brownian bridge.
Using well-known results by Pitman, Tanaka, Vervaat and Williams we are able to show that the bisection method has both a small
error and a small probability of failure.}

\section{Introduction and main result}

Imagine that we have a method of calculating values of a given function, $f:S \to \R$. 
But,  for each $x \in S$, our 
method for calculating $f(x)$ is computationally expensive.
The goal is to find $\min_{x \in S} f(x)$.
A good algorithm, balances the cost of further function evaluations.

We illustrate with a particular ``solvable'' example, which still has some typical complications.
We consider the example of finding the minimum of a Brownian bridge.

\subsection{Online minimization  in a Brownian bridge}

Let $B_{\mathbb{T}}$ denote the Brownian bridge, on $[0,1]$: $B_{\mathbb{T}}(0)=B_{\mathbb{T}}(1)=0$. Let $u \in (0,1)$ be such that
$$
	B_{\mathbb{T}}(u)\, =\, \min(\{B_{\mathbb{T}}(t)\, :\, t \in [0,1]\})\, .
$$
We seek an approximation of $u$ with an error of at most $\delta>0$, given.
In the following,
each ${\mathsf{StandardNormalRandomVariable}}$ is a new instance of a standard normal random variable.

\noindent
\underline{\bf Pseudo-code 1 : Brownian bridge fill-in function}\\
\underline{\sc Inputs:}\qquad $n \in \{1,2,\dots\}$ and $x_0,\dots,x_n \in \R$.\\
\underline{\sc Outputs:}\qquad $x_0',\dots,x_{2n}'$ equals $\mathsf{BBFI}(n,x_0,\dots,x_n)$.\\
\underline{\sc First for loop:}\qquad Repeat the following for $k \in \{0,\dots,n\}$.\\
\indent Let $x_{2k}' = x_k$.\\
\underline{\sc Second for loop:}\qquad Repeat the following for $k \in \{1,\dots,n\}$.\\
\indent Let $\displaystyle x_{2k-1}' = \frac{1}{2}\, x_{2k-2} + \frac{1}{2}\, x_{2k} + \frac{1}{\sqrt{2n}}\, \mathsf{StandardNormalRandomVariable}$.
\begin{gather*}
\hline
\end{gather*}

\vspace{-0.75cm}

\noindent
\underline{\bf Pseudo-code 2 : Initialization function}\\
\underline{\sc Inputs:}\qquad $d \in \{1,2,\dots\}$.\\
\underline{\sc Outputs:}\qquad $x_0,\dots,x_{2^d},t_0,\dots,t_{2^d}$ equals $\mathsf{Init}(d)$.\\
\underline{\sc Initialization:}\qquad Let $x_0=0$ and $x_1=0$.\\
\underline{\sc First for loop:}\qquad Repeat the following for $k \in \{0,\dots,2^d\}$.\\
\indent Let $t_{k} = k/2^d$.\\
\underline{\sc Second for loop:}\qquad Repeat the following for $r \in \{1,\dots,d\}$.\\
\indent Let $(x_0',\dots,x_{2^r}') = \mathsf{BBFI}(r,x_0,\dots,x_{2^{r-1}})$.\\
\indent Let $(x_0,\dots,x_{2^r}) = (x_0',\dots,x_{2^r}')$.

\newpage 

\noindent
\underline{\bf Pseudo-code 3 : Basic algorithm}\\
\underline{\sc Inputs:}\qquad $d$  for $N$, $d \in \{1,2,\dots\}$ and $d\geq 3$.\\
\underline{\sc Outputs:}\qquad A sequences of numbers $t^*_0,\dots,t^*_N$ and $\mathsf{Certificate}$.\\
\underline{\sc Initial steps:}\qquad Let $(t_0,\dots,t_{2^d},B_{\mathbb{T}}(t_0),\dots,B_{\mathbb{T}}(t_{2^d})) = \mathsf{Init}(d)$.\\
Find the $K(0)$ in $\{0,1,\dots,2^d\}$ such that 
$$
	B_{\mathbb{T}}(t_{K(0)})\, =\, \min(\{B_{\mathbb{T}}(t_k)\, :\, k \in \{0,1,\dots,2^d\}\})\, .
$$
Let $t^*_0 = t_{K(0)}$. Let $\mathsf{Certificate} =$ green-$\checkmark$.\\
For $k \in \{0,\dots,2^d\}$, let $\widetilde{B}^{(1)}(t_{k}) = 
B_{\mathbb{T}}(t_k+t_0^*-(1/2)\ \mathrm{mod}(1))-B_{\mathbb{T}}(t_0^*)$.\\
For $k \in \{0,\dots,2^{d-1}\}$ let $\widehat{B}^{(1)}(t_{k}) = \sqrt{2}\, \widetilde{B}^{(1)}(t_{2^{d-2}+k})$.\\
\underline{\sc For loop:}\qquad Repeat the following for $n \in \{1,2,\dots,N-1\}$.\\
\indent Let $({B}^{(n)}(t_0),\dots,{B}^{(n)}(t_{2^d}))$ equal $\mathsf{BBFI}(2^{d-1},\widehat{B}^{(n)}(t_0),\dots,\widehat{B}^{(n)}(t_{2^{d-1}}))$.\\
\indent Find the ${K}(n)$ in $\{0,1,\dots,2^d\}$ such that
$$
	{B}^{(n)}(t_{{K}(n)})\, =\, \min(\{{B}^{(n)}(t_k)\, :\, k \in \{0,1,\dots,2^d\}\})\, .
$$
\indent If ${K}(n)\not\in\{2^{d-1}-2^{d-3},\dots,2^{d-1}+2^{d-3}\}$, let  $\mathsf{Certificate} =$ red-X and {\bf abort}.\\
\indent Let $t^*_n = t_{K(n)}$.\\
\indent For $k \in \{0,\dots,2^{d-1}\}$, let 
$$
\widehat{B}^{(n+1)}(t_{k}) = \sqrt{2}\, \left({B}^{(n)}(t_{k+K(n)-2^{d-2}}) - {B}^{(n)}(t_n^*)\right)\, .\\
$$
\begin{gather*}
\hline
\end{gather*}

\vspace{-0.75cm}

\begin{figure}
\begin{center}
\begin{tikzpicture}
\draw (0,0) node[] {\includegraphics[width=12cm]{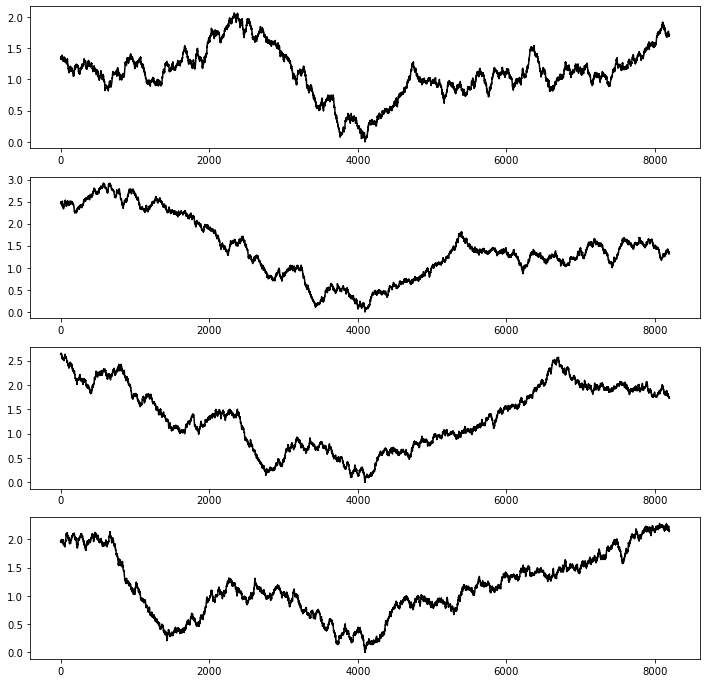}};
\end{tikzpicture}
\caption{Output of the basic algorithm. We chose $d=14$ and $N=4$.
For each $n \in \{1,\dots,4\}$, we have plotted
$\widehat{B}^{(n)}(t_k)$ against $k$, for  $k \in \{0,\dots,2^{d-1}\}$. This is a random sample path.
\label{fig:1}}
\end{center}
\end{figure}

See Figure \ref{fig:1} for a sample of the main algorithm.
\begin{theorem} Let $U_0 = t_0^*$. For each $n \in \{1,\dots,N\}$, let  $U_n =U_{n-1}+ \left(t^*_n -\frac{1}{2}\right) 2^{-n}$. Then both
 $\mathbf{P}(\mathsf{Certificate}=\text{red-X})$ and  $\mathbf{P}(\operatorname{dist}(u,U_N+\Z)>2^{-N})$ are  $O(N\, d^{1/2}/2^{d/2})$
as $N,d \to \infty$.
\end{theorem}


The run-time is $T = O(N2^d)$. 
By the theorem,  if we want to bound the probability of failure to be at most $p_0$, we must take $(2^d/d) = O(N^2/p_0^2)$:
i.e., $2^d = O((N/p_0)^2 \ln(N/p_0))$.
Therefore, since $\delta = 2^{-N}$, that means $N = O(|\ln(\delta)|)$. So we have the run time $T = O(|\ln(\delta)|^3 p_0^{-2} \ln(|\ln(\delta)|p_0^{-1}))$.

We recall a well-known result for the analogous discrete problem: determining the minimum number of look-up steps 
to find the minimum for a simple random walk $X_0,X_1,\dots,X_n$. This was considered in \cite{Odlyzko} (and later
 \cite{Chassaing}). The result is that the number of look-up steps is 
$(c+o(1)) \sqrt{n}$ for $c = \sqrt{8/\pi}\, \ln(2)$. The two problems are different.

We review the following well-known facts, which we use in the sequel.
\begin{itemize}
\item[$\bullet$]{\bf Fact1:}  For Brownian motion, if we have two times $t_1<t_2$ and we condition on $B(t_1)=x$ and $B(t_2)=y$, then the conditional 
distribution of $\min(\{B(t)\, :\, t_1<t<t_2\})$ is known by the reflection principle.
\item[$\bullet$]{\bf Fact2:} Vervaat proved the process-wise equivalence between the Brownian excursion process and a transformation of the Brownian
bridge, where for $(t,B_t) \in  (\R/\Z)\times \R$, you ``shift'' the minimum of the Brownian bridge to the origin $(0,0)$ \cite{Vervaat}.
\item[$\bullet$]{\bf Fact3:} Williams proved that the Brownian excursion process (BEP) is a Bessel(3) bridge \cite{Williams}.
\item[$\bullet$]{\bf Fact4:} On all of $[0,\infty)$, Pitman proved the distributional equivalence between the past (running) maximum of Brownian motion started at $(0,0)$
and the future (tail) minimum of the Bessel(3) process started at $(0,0)$ \cite{Pitman}
\end{itemize}

\section{A second certificate}

\label{sec:certificate}

The following certificate is not essential.
The basic algorithm above is all that is necessary. But it helps in the proof.

Let us denote by $\widetilde{B} : \R \to \R$ the function $\widetilde{B}(t) = B_{\mathbb{T}}(t\ \mathrm{mod}(1))$.
For each $n \in \{1,\dots,N\}$, let $\overline{B}^{(n)} : \R \to \R$, be the function such that
$$
\overline{B}^{(n)}(t)\, =\, 2^{(n-1)/2} \left(\widetilde{B}\left(\frac{t}{2^{n-1}}+\sum_{r=0}^{n-1} t^*_k 2^{-k}-\frac{1}{2^n}\right)
	-  \widetilde{B}\left(\sum_{r=0}^{n-1} t_k^* 2^{-k}\right)\right)\, .
$$
As long as the certificate is a green-$\checkmark$ at the $n$th step,
this is just the extension of $B^{(n)}$ to continuous values: 
$\forall k \in \{0,1,\dots,2^d\}$ we have $\overline{B}^{(n)}(t_k) = B^{(n)}(t_k)$.

At the risk of pedantry, let $\mathcal{F}_{r}$ be the $\sigma$-algebra generated by $\widetilde{B}(k/2^r)$ for all $k \in \{0,1,\dots,2^r\}$.
Then the random variables $t_0^*,\dots,t_{n-1}^*$ are $\mathcal{F}_{d+n-1}$-measurable.
Then by the usual Markov/Gibbs property of Brownian motion, for the portion of the path $(\overline{B}^{(n)}(t))_{t \in (t_{k-1},t_k)}$, if we condition on $\mathcal{F}_{d+n-1}$,
the conditional distribution is the same as if we just condition on $B^{(n)}(t_{k-1})$ and $B^{(n)}(t_{k})$.
Moreover, by the reflection principle, Fact 1,
\begin{equation*}
\begin{split}
&
	\mathbf{P}\left(\min\left(\left\{\overline{B}^{(n)}(t) : t \in (t_{k-1},t_{k})\right\}\right) \leq z\, \Big|\, 
B^{(n)}(t_{k-1})=x,B^{(n)}(t_{k})=y\right)\\ 
&\hspace{2cm}
=\, \begin{cases}
1\, , & \text{if $z \geq \min(\{x,y\})$,}\\
\exp(-2^{d+1}(z-x)(z-y))\, , & \text{ otherwise.}
\end{cases}
\end{split}
\end{equation*}
We may invert this formula: this equals $p \in (0,1)$ if and only if 
$$
	z\, =\, \frac{x+y}{2} - \sqrt{\frac{(x-y)^2}{4} - \frac{\ln(p)}{2^{d+1}}}\, .
$$
This gives a method for simulating the minimum, by a well-known procedure.

For the following pseudo-code, each instance of $\mathsf{Uniform}(0,1)$ is a new, independent copy of a uniform pseudo-random variable
in $(0,1)$.
The pseudo-code describes a sub-routine which is supposed to be added-on to the previous pseudo-code at the end of each step of the main $n$-For-loop.
We assume that Certificate 1 was a green-$\checkmark$. Otherwise, we should have already had an early-abort.

\bigskip

\noindent
\underline{\bf Pseudo-code 2 : Add-on algorithm at the end of the $n$th step}\\
\underline{\sc Inputs:}\qquad ${B}^{(n)}(t_k)$ for $k \in \{0,\dots,2^d\}$ and $K(n)$.\\
\underline{\sc Outputs:}\qquad $\mathsf{Certificate}$-2.\\
\underline{\sc Initial step:}\qquad 
Let Certificate 2 $=$ green-$\checkmark$.\\
\underline{\sc For loop:}\qquad Repeat the following for $k \in \{1,\dots,2^d\}$.\\
If $\{t_{k-1},t_k\} \not\subset \{2^{d-2},\dots,2^{d-2}+2^{d-1}\}$\\
\indent then let $x = B^{(n)}(t_{k-1})$, $y=B^{(n)}(t_k)$, and\\
\indent let
$$
	m\, =\, \frac{x+y}{2} - \sqrt{\frac{(x-y)^2}{4} - \frac{\ln(\mathsf{Uniform}(0,1))}{2^{d+1}}}\, .
$$
\indent If $m\leq B^{(n)}(t_{K(n)})$, change Certificate 2 to a red-X and {\bf abort}.

\begin{gather*}
\hline
\end{gather*}

\vspace{-0.75cm}

\begin{figure}
\begin{center}
\begin{tikzpicture}
\draw (0,0) node[] {\includegraphics[width=12cm]{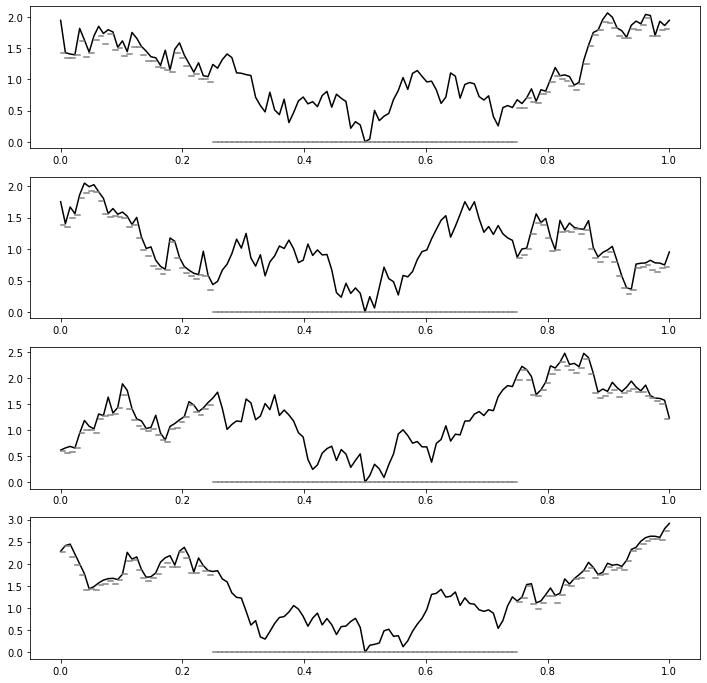}};
\end{tikzpicture}
\caption{Output of the basic algorithm along with a graphical depiction of the $m$-values entering the 2nd certificate. We chose $d=7$ and $N=4$.
(We used a lower $d$ to better demonstrate the $m$-values.)
For each $n \in \{1,\dots,4\}$, we have plotted
$\sqrt{2}({B}^{(n)}(t_k)-B^{(n)}(t_{K(n)}))$ in black, as well as the $m$-values for the 2nd certificate in gray. We have $m$-values in the first and fourth subintervals. In the middle two subintervals, we just use $0$. Note that the restriction of the black plot to the middle 2 subintervals
gives $\widehat{B}^{(n)}(t_k)$ for $k \in \{0,\dots,2^{d-1}\}$.
If we ever had the $m$-values of the first and fourth subintervals being less than 0, the 2nd certificate would return ``red X'' and abort.
\label{fig:2}}
\end{center}
\end{figure}

See Figure \ref{fig:2} for a graphical representation including an explantion of the event which would trigger Certificate-2 to be a red-X.

\section{Proofs}
From 
Vervaat's theorem, Fact 2, we know that $(\widetilde{B}(u+t)-\widetilde{B}(u))_{t \in [0,1]}$ is a Brownian excursion process (BEP)
independent of $u$.
In turn, the BEP is a Bessel bridge by Williams's theorem, Fact 3. This means that if $Y : [0,\infty) \to [0,\infty)$ is a standard Bessel(3) process,
started at $(0,0)$, then the BEP is equally distributed to
$$
	(1-t)\, Y\left(\frac{t}{1-t}\right)\ \text{ for }\ 0<t<1\, .
$$
Let $u_n$ be the arg-min of $\overline{B}^{(n)}$ on $[0,2^{n-1}]$. Then we know that the process
$(\overline{B}^{(n)}(u_n+t)-\overline{B}^{(n)}(u))_{t \in [0,2^{n-1}]}$
is equally distributed to
$$
	(1-t2^{-n+1})\, Y^{(n)}\left(\frac{t}{1-t2^{-n+1}}\right)\ \text{ for }\ 0<t<2^{n-1}\, ,
$$
where $Y^{(n)}$ is $2^{(n-1)/2} Y(2^{-n+1}\, \cdot)$. Then by scaling invariance
of Bessel-3, the process $Y^{(n)}$ is distributed identically to the process $Y$.

We note that the identification of $Y$ as the Euclidean norm of 3d Brownian motion shows that, for any $z>0$, we have
$$
	\mathbf{P}(\max(\{Y(t)\, :\, 0\leq t\leq 1\})\geq z)\, \leq\, 6 \int_{z/3}^{\infty} \frac{e^{-x^2/2}}{\sqrt{2\pi}}\, dx\, \leq\, 
	\frac{18 e^{-z^2/18}}{z\, \sqrt{2\pi}}\, .
$$
See, for instance, Problem 1 from \cite{McKean}.
The opposite bound follows from Pitman's theorem, Fact 4.
\begin{lemma}
For any $z>0$, we have
$$
\mathbf{P}(\min(\{Y(t)\, :\, t\geq 1\})\leq z)\, =\, 1-2\int_{z}^{\infty}\, \frac{e^{-x^2/2}}{\sqrt{2\pi}}\, dx\, 
\leq\, \frac{2z}{\sqrt{2\pi}}\, .
$$
\end{lemma}
By scaling, $Y$ is equidistributed to $c^{-1/2} Y(c\, \cdot)$. So, for example
$$
\mathbf{P}(\max(\{Y(t)\, :\, 0\leq t\leq \varepsilon\})\geq z )\, =\, 
\mathbf{P}(\max(\{Y(t)\, :\, 0\leq t\leq 1\})\geq z\varepsilon^{-1/2})\, .
$$
Then by the two inequalities and the subset bound we have the following.
\begin{corollary} For any $\varepsilon>0$, $a \in (\varepsilon,\infty)$ and $\lambda \in (0,\infty)$, let us define
$$
	\mathcal{P}(\varepsilon,a,\lambda)\, =\, 
	\mathbf{P}\Big(\min\big(\{Y(t)\, :\, t\geq a\}\big) \leq \lambda^{1/2} \max\big(\{Y(t)\, :\, 0\leq t\leq \varepsilon\}\big)\Big)\, .
$$
Then 
$$
	\mathcal{P}(\varepsilon,a,\lambda)\, \leq\, \frac{1}{3}\, \sqrt{\frac{2\varepsilon\lambda}{\pi a}\, \ln\left(\frac{a}{\varepsilon\lambda}\right)}\, .
$$
\end{corollary}
Let $\mathcal{P}(\varepsilon,a) = \mathcal{P}(\varepsilon,a,1)$ so that
$$
	\mathcal{P}(\varepsilon,a)\, \leq\, \frac{1}{3}\, \sqrt{\frac{2\varepsilon}{\pi a}\, \ln\left(\frac{a}{\varepsilon}\right)}\, .
$$
As long as the aspect ratio is eccentric, the probability is low. 

Consider $\widetilde{B}$. For some $\Delta \in \{1,2,\dots\}$, let us choose the lattice scaling to be $2^{-\Delta}$.
Let $t_k(\Delta) = k/2^\Delta$ for $k \in \Z$. Now, let $\kappa_n(\Delta)$ be that number $k \in \{0,\dots,2^{\Delta}\}$
such that $|t_{k}(\Delta)-u_n|$ is minimal. In particular,
$$
	|t_{\kappa_n(\Delta)}-u_n|\, \leq\, \frac{1}{2^{\Delta}}\, .
$$
So we have
\begin{equation*}
\begin{split}
	&\mathbf{P}\Big(\exists k \in \{0,\dots,2^{\Delta}\}\, \text{ s.t. } (|t_k-t_{\kappa_n(\Delta)}|>a) \wedge 
\big(\overline{B}^{(n)}(t_k)<\lambda^{1/2} \overline{B}^{(n)}(t_{\kappa}(\Delta))\big)\Big)\\
&\hspace{7.5cm} \leq\, 2\cdot \mathcal{P}(2^{-d},a/2,4\lambda)\, ,
\end{split}
\end{equation*}
where we take into account the fact that the multiplier $(1-t2^{-n+1})$ could be as small as $1/2$ (so $a$ changes to $a/2$) and the time
stretching $t\mapsto t/(1-t2^{-n+1})$ could scale $t$ by as much as $2$ (so $\lambda$ changes to $4\lambda$), as long as $0<t<2^{n-2}$ is in the Voronoi neighborhood
of the left endpoint of $[0,2^{n-1}]$ (and the
Voronoi neighborhood of the right endpoint satisfies the same inequalities by symmetry).

Now we run an induction argument, using the subset inequality.
First, suppose that all the Certificate-1's and Certificate-2's give green-$\checkmark$s for levels up to but not including $n$
(for some $n \in \{0,\dots, N-1\}$).
Also, suppose that Certificate-1 at level $n$ is a green-$\checkmark$.
That all is going to operate as our induction hypothesis.

Now suppose that we have $|u_n-1/2|>7/32$.
Then this would mean $|u_{n-1}-1/2|>7/16$.
But we had $|t^*_{n-1}-1/2|<3/8=6/16$ because of the green-$\checkmark$ for Certificate-1 at the $(n-1)$st step.
So we have an event whose probability is bounded by $2\cdot\mathcal{P}(2^{-d+1},c/16,C)$, for $C\leq 4$ and $c\geq 1/2$.
This is $O(\sqrt{d2^{-d}})$.

Then, on the complement of this event, we have that $|u_n-1/2|<7/32$.
So the probability that $\min(\{\overline{B}^{(n)}(t):|t-1/2|\geq 1/4\})<\overline{B}^{(n)}(t_{K(n)})$ is
an event whose probability is at most $2\cdot\mathcal{P}(2^{-d},c/32,C)$, which is $O(\sqrt{d2^{-d}})$.
(Recall $B^{(n)}(t_{K(n)})$ is even smaller than or equal to $B^{(n)}(t_{\kappa_n(d)})$.)
But if this event does not happen, we have Certificate-2 is a green-$\checkmark$.

Now, let us get the bound on the complementary probability for Certificate-1 at level $(n+1)$ being a green-$\checkmark$.
Let us assume the opposite, that $|t_{K(n+1)}-1/2|>1/8$ in order to show the probability is small.
We note that, by construction $\overline{B}^{(n+1)}(t_{2k})\geq \overline{B}^{(n+1)}(1/2)$ for all $k \in \{0,\dots,2^{d-1}\}$.
So, if we have $|u_{n+1}-1/2|>1/16$ then that is an event whose probability is bounded by 
$2\cdot \mathcal{P}(2^{-d+1},c/16,C)$,
which is $O(\sqrt{d2^{-d}})$.
But if $|u_{n+1}-1/2|<1/16$, then that means the event that $|t_{K(n+1)}-1/2|>1/8$ has probability
bounded by $2\cdot \mathcal{P}(2^{-d},c/16,C)$, which is $O(\sqrt{d2^{-d}})$.
This completes the induction step.

Since we have Certificate-1 at level-0 is a green-$\checkmark$ by default, the initial step is trivial.
So, by induction, we have proved that all the certificates are green-$\checkmark$'s for $n$ up to and including $N-1$
with a total complementary probability at most $O(N\, \sqrt{d2^{-d}})$.
But if all the certificates are green-$\checkmark$'s, for Certificates-1 and Certificates-2, 
then $\operatorname{dist}(u,U_N+\Z)$ is at most $2^{-N}$.

\section{Summary and outlook}

We have outlined a computer program for finding the minimum of Brownian bridge online. 
By online, we mean that the search algorithm operates simultaneously with the recursive fractal
construction of the Brownian bridge.
Python codes that implement this,  and in particular which generated the images in Figures \ref{fig:1} and \ref{fig:2}
are available freely, online, at GitHub \cite{StarrRepository}.

We are led to an alternative description of the Bessel-3 process.
Consider a spatial process $X : \R \to \R$. For example, this could be a Bessel-3 process.

Now, let $T_0$ be a spatial Poisson process on $\R$ of rate $1$, conditioned on $0 \in T_0$.
(in other words, let $S_0$ be a spatial Poisson process on $\R$ of rate 1, and let $T_0 = S_0 \cup \{0\}$ by
the memoryless property.)
Enumerate the points of $T_0$ as $\{\dots,t_{-1},t_0,t_1,\dots\}$ where $t_0=0$.
Define a spatial process $\mathcal{X} : \R \to \R$ as follows.
Let $\mathcal{X}(t_n) = X(t_n)$ for each $n \in \Z$.
For each interval $(t_{n-1},t_n)$, let $\mathcal{X}$ on $(t_{n-1},t_n)$ be Brownian motion conditioned on $\mathcal{X}(t_{n-1})=X(t_{n-1})$
and $\mathcal{X}(t_n)=X(t_n)$.

Now let us use a construction, which is sometimes useful \cite{StarrVermesiWei}.
Let $\mathcal{T}$ be a space-time Poisson process on $\R \times [0,\infty)$ with area rate 1.
Let $T_t = \{x\, :\, (\exists s \in [0,t)\ \text{s.t. } (s,xe^{-s}) \in \mathcal{T}) \vee (xe^{-t} \in T_0)\}$.
Then $T_t$ is a spatial Poisson process on $\R$ of rate $1$ for all $t\geq 0$. But, more importantly, the $T_t$'s fit together
in a good way (as in the defunct steady state cosmological model \cite{Weinberg}).
The process $T_t$ is expanding at a uniform rate, and new particles are being injected according to a space-time Poisson process
in order to maintain constant average density.

Now, consider a function $Y_t : T_t \to \R$ defined as follows:
$$
\forall x \in T_0\, ,\ Y_t(e^tx) = e^{t/2} X(x)\ ,
$$
and
$$
\forall (x,s) \in \R \times (0,t]\, ,\ Y(e^{t-s}x) = e^{t/2} \mathcal{X}(e^{-s}x)\, .
$$
We are not done, yet.
There is an extra part of the Markovian dynamics.
For each $t \in [0,\infty)$, let $x^*(t) \in T_t$ be the arg-min of $Y_t$.
Then let $\widetilde{T}_t = \{x-x^*(t)\, :\, x \in T_t\}$ and let $\widetilde{Y}_t : \widetilde{T}_t \to \R$ be given by
$$
	\widetilde{Y}_t(y)\, =\, Y_t(y+x^*(t)) - Y_t(x^*(t))\, .
$$
Now let $\mathcal{G}_t = \{(y,\widetilde{Y}_t(y))\, :\, y \in \widetilde{T}_t\}$ be the graph.

We think an interesting conjecture is that $\mathcal{G}_t$ is (well-defined and) stationary if and only if $X$ is a Bessel-3 process.
If that were proved, then one might investigate the convergence rate.

\section*{Acknowledgments}
S.S.\ was funded by a Simons collaboration grant, which benefitted this project.
S.S.\ wishes to thank Arvind Singh for a useful discussion, in particular for suggesting to look at \cite{Pitman}.

\label{sec:Outlook}

\baselineskip=12pt
\bibliographystyle{plain}

\end{document}